\documentclass[reqno]{amsart}

\usepackage{ mathtools, amssymb, mathdots, bbm, mleftright, bm, subcaption, graphicx}
\usepackage[dvipsnames]{xcolor}
\definecolor{mygreen}{HTML}{B01C2E}

\usepackage[unicode]{hyperref}

\usepackage{ tikz, tikz-cd }
\usetikzlibrary{shapes.misc, calc, patterns}

\usepackage[shortlabels]{ enumitem }
\usepackage{graphicx}

\usepackage{ amsthm, thmtools, hyperref, url }
\usepackage[nameinlink]{ cleveref }

\numberwithin{equation}{section}
\newtheorem{theorem}{Theorem}[section]
\newtheorem*{theorem*}{Theorem}
\newtheorem{corollary}[theorem]{Corollary}

\newtheorem{proposition}[theorem]{Proposition}

\newtheorem{lemma}[theorem]{Lemma}

\theoremstyle{definition}
\newtheorem{definition}[theorem]{Definition}

\newtheorem{remark}[theorem]{Remark}
\newtheorem{example}[theorem]{Example}

\newcommand{\Z}{\mathbb{Z}}
\newcommand{\Q}{\mathbb{Q}}

\newcommand{\F}{\mathbb{F}}
\newcommand{\N}{\mathbb{N}}

\DeclareMathOperator{\Sym}{Sym}
\DeclareMathOperator{\SL}{SL}

\renewcommand{\leq}{\leqslant}

\renewcommand{\geq}{\geqslant}

\newcommand{\qbinom}[2]{\genfrac{[}{]}{0pt}{}{#1}{#2}}
\newcommand{\Qbinom}[2]{\genfrac{\{}{\}}{0pt}{}{#1}{#2}}

\newcommand{\B}{\mathcal{B}}
\newcommand{\A}{\mathcal{A}}
\renewcommand{\sl}{\mathfrak{sl}}

\def\C{\mathbb{C}} 
\def\U{\mathcal{U}}

\newcounter{thmlistcnt}
\newenvironment{thmlist}%
	{\setcounter{thmlistcnt}{0}%
	\begin{list}{\emph{(\roman{thmlistcnt})}}{%
		\usecounter{thmlistcnt}%
		\setlength{\topsep}{0pt}%
		\setlength{\leftmargin}{0pt}%
		\setlength{\itemsep}{-3pt}%
		\setlength{\labelwidth}{17pt}
		\setlength{\itemindent}{30pt}}%
	}%
	{\end{list}}%

\newcounter{defnlistcnt}
\newenvironment{defnlist}%
	{\setcounter{defnlistcnt}{0}%
	\begin{list}{(\alph{defnlistcnt})}{%
		\usecounter{defnlistcnt}%
		\setlength{\topsep}{-3pt}%
		\setlength{\leftmargin}{0pt}%
		\setlength{\itemsep}{0pt}%
		\setlength{\labelwidth}{25pt}
		\setlength{\itemindent}{30pt}}%
	}%
	{\end{list}}%

		\title[A new bijective proof of the $q$-Pfaff--Saalsch\"utz identity]{A new bijective proof \\ of the $q$-Pfaff--Saalsch\"utz identity \\ with applications to quantum groups} 
		
		\author[Á.~Gutiérrez]{\'Alvaro Guti\'errez}
		\author[Á.~Martínez]{\'Alvaro L.~Mart\'inez}
		\author[M.~Szwej]{Micha{\l} Szwej}
		\author[M.~Wildon]{Mark Wildon}
		
		\address{
			University of Bristol,	Bristol, UK\\
		}
        \email[ÁG]{a.gutierrezcaceres@bristol.ac.uk}
        \email[MS]{michal.szwej@bristol.ac.uk}
        \email[MW]{mark.wildon@bristol.ac.uk}
		\address{
			Columbia University, New York, USA
		}
        \email[ÁM]{alm2297@columbia.edu}

\begin{document}

		\begin{abstract}
			We present a combinatorial proof of the $q$-Pfaff--Saalschütz identity by a composition of explicit bijections, in which $q$-binomial coefficients are interpreted as counting subspaces of $\mathbb{F}_q$-vector spaces. As a corollary, we obtain a new multiplication rule for quantum 
			binomial coefficients and hence a new presentation of 
			Lusztig's integral form $\mathcal{U}_{\mathbb{Z}[q, q^{-1}]}(\mathfrak{sl}_2)$ of 
			the Cartan subalgebra of the quantum group $\mathcal{U}_q(\mathfrak{sl}_2)$.\smallskip


		\end{abstract}

	    \maketitle
	
		
		
		\section{Introduction}
		The study of binomial identities of the form $\binom{\cdot}{\cdot}\binom{\cdot}{\cdot}=\sum\binom{\cdot}{\cdot}\binom{\cdot}{\cdot}\binom{\cdot}{\cdot}$ dates back to an identity stated by Shih-Chieh Chu \cite{Tak73} in 1303. Multiple identities of this form were discovered around the 1970s by Nanjundiah \cite{Nan58}, Stanley \cite{Sta70}, Bizley \cite{Biz70}, Tak\'acs \cite{Tak73} and Sz\'ekely \cite{Sze85}. In the setting of hypergeometric series \cite[\S5.5]{GKP89}, it becomes clear that each of them is an instance of the Pfaff--Saalsch\"utz identity \cite{Zen89} --- in other words, they are all equivalent 
		by a simultaneous linear change of variables: see Examples~\ref{ex:Stanley} and~\ref{ex:Zeil},
		and Remark~\ref{remark:hypergeometric}. Each of these identities lifts to a $q$-analogue of the form $\qbinom{\cdot}{\cdot}_q\qbinom{\cdot}{\cdot}_q=\sum q^{\,\cdot}\qbinom{\cdot}{\cdot}_q\qbinom{\cdot}{\cdot}_q\qbinom{\cdot}{\cdot}_q$, with the same equivalences holding. We may therefore refer to any identity of this form as a
		\emph{$q$-Pfaff--Saalsch\"utz identity}. The combinatorial proofs to date of  $q$-Pfaff--Saalsch\"utz identities are due to Andrews and Bressoud \cite{And84}, Goulden \cite{Gou85}, Zeilberger \cite{Zei87}, Yee \cite{Yee08}, and Schlosser and Yoo \cite{SY17}.

		Let $\qbinom{n}{k}_q$ denote the $q$-binomial coefficient, defined in \S\ref{sec:background} for
		$n$, $k \in \N_0$ as a polynomial in $q$.
		In \S\ref{sec:new-bijective-proof} we present a new combinatorial proof of the $q$-Pfaff--Saalsch\"utz identity in the following form.
		
		\begin{theorem}[$q$-Pfaff--Saalsch\"utz identity]\label{thm:q-nan}
			Let $m,s,t\in\mathbb{N}_0$. If $e \in \Z$ and \break$-t \le e \le s$ then
			\begin{align*}
				\qbinom{m}{t}_q\qbinom{m+e}{s}_q=\sum_{j\geq 0}q^{(s-j)(t+e-j)}\qbinom{t+e}{j}_q\qbinom{s-e}{s-j}_q\qbinom{m+j}{s+t}_q.
			\end{align*}
			
		\end{theorem}

		With the exception of Zeilberger's proof \cite{Zei87}, 
		which relies on an earlier bijection due to Foata \cite[page 233]{Foa65}, our 
		proof is shorter than all those in the literature to date.
		Moreover, our proof uses a combinatorial interpretation of $q$-binomial coefficients not previous used in this context:
		if $q$ is a prime power then the
		 $q$-binomial coefficient $\qbinom{n}{k}_q$ is the 
		number of $k$-dimensional subspaces of~$\mathbb{F}_q^n$ (see \S\ref{sec:background}).
		Finally we emphasise that while our proof has multiple steps, each step is bijective.

		Our second main result lifts Theorem~\ref{thm:q-nan} to an identity in the quantum group $\U_q(\sl_2)$. In \cite{Lus88}, Lusztig defines an integral form $\U_\A(\sl_2)$ for $\U_q(\sl_2)$ over $\A=\Z[q,q^{-1}]$. It can be thought of as a deformation of Kostant's $\mathbb{Z}$-form for $\U(\sl_2)$, and it admits a triangular decomposition
		\[\U_\A(\sl_2) = \U_\A^-(\sl_2)\otimes_{\A}\U_\A^0(\sl_2)\otimes_{\A}\U_\A^+(\sl_2).\]
		To construct the subalgebra $\U_\A^0(\sl_2)$, Lusztig uses the elements
		\begin{align}\label{formula:q-kostant}
			\qbinom{K;c}{t}=\frac{[K;c][K;c-1]\cdots [K;c-t+1]}{\{t\}_q!},
		\end{align}
		for $c \in \Z$ and $t \in \N_0$, where
		$[K;a]= \frac{q^aK-q^{-a}K^{-1}}{q-q^{-1}}$ and $\{t\}_q!$ is the quantum factorial defined in  \S\ref{sec:q-to-quantum}. A lift of the identity in Theorem~\ref{thm:q-nan} to the quantum group gives the multiplication rule for these elements. 
		In the following theorem, $\Qbinom{n}{k}_q$ denotes the quantum binomial coefficient, defined as a Laurent
		polynomial in $q$ in \S\ref{sec:q-to-quantum}.
		
		\begin{theorem}\label{thm:multiplication_rule}
			Let $b,c\in\Z$ and $s,t\in\N_0$. If $t-c+b\geq0$ and $s-b+c\geq 0$, then the elements defined in~\eqref{formula:q-kostant} satisfy the following multiplication rule:
			\begin{align*}
				\qbinom{K; c}{t}\qbinom{K; b}{s}=\sum_{i\ge 0} \Qbinom{t-c+b}{i-c}_q\Qbinom{s-b+c}{i-b}_q\qbinom{K;i}{t+s}.
			\end{align*}
			
		\end{theorem}
		
		Lusztig shows that the elements $K^\delta\qbinom{K; 0}{t}$ for $t\ge 0$ and $\delta\in \{0,1\}$ form an $\A$-basis for $\U_\A^0(\sl_2)$. 
		In Proposition~\ref{prop:binomialbasis}, we show that the elements $\qbinom{K;\hskip1pt c}{t}$ for $t\ge 0$ and $c\in \{0,1\}$ also form an $\A$-basis for $\U_\A(\sl_2)$. This gives our new description of the multiplication in $U^0_\A(\sl_2)$ in Theorem~\ref{thm:cartanpresentation}. As a corollary, we obtain a presentation of $U^0_\A(\mathfrak{g})$ for an arbitrary Kac--Moody algebra $\mathfrak{g}$ of finite rank: see Corollary~\ref{cor:allranks}.
		
		\subsubsection*{Outline}
		In \S\ref{sec:background} we give foundational results about $q$-binomial coefficients, including bijective proofs of certain basic identities using the vector space interpretation. In~\S\ref{sec:new-bijective-proof} we present our new combinatorial proof of Theorem~\ref{thm:q-nan}. In \S\ref{sec:q-to-quantum} we give the 
		quantum version of this identity, later used to prove Theorem~\ref{thm:multiplication_rule}. 
		The quantum group~$\U_q(\sl_2)$ and Lusztig's integral form are briefly reviewed
		in \S\ref{sec:quantum}. Finally in \S\ref{sec:quantum-new} we prove Theorem~\ref{thm:multiplication_rule} and, 
		as a corollary, we obtain a new presentation of Lusztig's integral~form.

		\section{\texorpdfstring{$q$}{q}-binomials and subspaces of \texorpdfstring{$\F_q^n$}{F\_{}q\^{}n}}\label{sec:background}
		
		In this section, we state foundational results on $q$-binomial coefficients. 
		We include certain bijective proofs to justify our claim that each step in the proof of Theorem~\ref{thm:q-nan} is bijective.
		We refer the reader to \cite{KC02} for further background.		
		
		We define the $q$\textit{-integer} by
		\begin{align*}
			[n]_q=q^{n-1}+\dots+q+1=\frac{q^n-1}{q-1}
		\end{align*}
		for $n \in \N_0$. 
		The definitions of $q$-factorial and $q$-binomial coefficient follow naturally: 
		$[n]_q!=[n]_q[n-1]_q\ldots[1]_q$ and
		\begin{align}\label{eq:defn}
			\qbinom{n}{k}_q=\frac{[n]_q!}{[k]_q![n-k]_q!}
		\end{align}
		for $n, k \in \N_0$.
		We set $\qbinom{n}{k}_q=0$ whenever $n<k$.
		
		\begin{remark}
			Under this convention, the sum in Theorem~\ref{thm:q-nan} is finite, with non-zero summands for $\max\{0,e\}\leq j\leq \min\{t+e, s\}.$
		\end{remark}
		
		The $q$-binomial coefficients are the generating functions
		for subsets of~$\N_0$
		by the sum of their elements, partitions in an $(n-k)\times k$ rectangle, and inversions in permutations: 
		see \cite[Chapter~1]{Sta12}. We refer the reader to \cite{Pra11} for further background including an introduction to positivity results.
		Instead we use the following combinatorial interpretation, which could in the context
		of the proof of Theorem~1.1 be taken as the
		\emph{definition} of $\qbinom{n}{k}_q$.
		For the remainder of this section, $q$ is a prime power. 
		
		\begin{proposition}\label{prop:q-binom-defn}
			The number of $k$-dimensional subspaces of $\mathbb{F}_q^n$ equals $\qbinom{n}{k}_q$.
		\end{proposition}
		
		\begin{proof}
			See \cite[Thm.~7.1]{KC02}.
		\end{proof}
		
%
		\begin{example}
			Let $n=2$, $k=1$. By definition, $\qbinom{2}{1}_q=\frac{q^2-1}{q-1}=q+1$. When $q=3$, there are  four lines through the origin in~$\mathbb{F}_3^2$, as shown in Figure \ref{fig:f3-example}. 
						More generally, there are $(q^n-1)/(q-1) = 1 + q + \cdots + q^{n-1} = [n]_q = \qbinom{n}{1}_q$ 
			lines through the origin in $\F_q^n$.
			\begin{figure}[h]
				\centering
				\begin{tikzpicture}
					\foreach \x in {0,1,2} {
						\foreach \y in {0,1,2} {
							\fill[black] (\x, \y) circle (2.5pt);
						}
					}
					\draw[loosely dotted, thick] (-.25,0) -- (2.25,0);
					\draw[loosely dashed, thick] (0,-.25) -- (0,2.25);
					\draw[loosely dashdotted, thick] (0,0) -- (1,2) -- (2,1) -- (0,0);
					\draw[thick] (-.25, -.25) -- (2.25, 2.25);
				\end{tikzpicture}
				\caption{The four one-dimensional subspaces of $\mathbb{F}_3^2$, visualizing~$\F_3^2$~as~a~subset~of~the~plane.}
				\label{fig:f3-example}
			\end{figure}
		\end{example}
		
		This interpretation allows to lift  bijective proofs of classical binomial identities, in which
		$\binom{n}{k}$ counts the number of $k$-subsets of $\{1,\ldots, n\}$, into bijective proofs of their $q$-analogues.
		A typical instance is the following lemma. While it has a one-line algebraic proof,
		we give a bijective proof in order that,
		when it is used in our proof of Theorem~\ref{thm:q-nan}, the entire
		proof is by a composition of bijections.
		Throughout we use $\leqslant$ for containment of subspaces.

		\begin{lemma}\label{lem:n-k_choose_l-k}
			Let $U$ be a $k$-dimensional subspace of $V=\F_q^n$. The number of \break$\ell$-dimensional subspaces $W$ of $V$, such that $U\le W\le V$, equals $\qbinom{n-k}{\ell-k}_q$.
		\end{lemma}
		\begin{proof}
			There is a natural bijection between the 
			$\ell$-dimensional subspaces $W$ of $V$ and the $(\ell-k)$-dimensional subspaces of the quotient space $V/U$,
			 given by
			$W\longleftrightarrow W/U$. By Proposition~\ref{prop:q-binom-defn} there are $\qbinom{n-k}{\ell-k}$ such $(\ell-k)$-dimensional subspaces.
		\end{proof}
		\noindent We deduce a $q$-analogue of the trinomial revision identity $\binom{n}{l}\binom{l}{k}=\binom{n}{k}\binom{n-k}{l-k}$.
		
		\begin{proposition}\label{lemma:subset_of_subset}
			Let $k\leq \ell\leq n$ be non-negative integers. Then
			\begin{align*}
				\qbinom{n}{\ell}_q\qbinom{\ell}{k}_q=\qbinom{n}{k}_q\qbinom{n-k}{\ell-k}_q.
			\end{align*}
		\end{proposition}
		
		\begin{proof}
			By Proposition~\ref{prop:q-binom-defn}, the left-hand side is the number of pairs of subspaces $(U,W)$ of $V=\mathbb{F}_q^n$ satisfying $U\leq W\leq V$, $\dim U=k$ and $\dim W=\ell$. By Lemma~\ref{lem:n-k_choose_l-k}, the same pairs are counted by the expression on the right-hand side, so the two are equal.
%
		\end{proof}
		
		Another $q$-analogue of a basic binomial identity, which again 
		we shall need in the proof of Theorem~\ref{thm:q-nan} is as follows. Again we give a bijective~proof.
		
		\begin{proposition}\label{lemma:n-k}
			Let $k\leq n$ be non-negative integers. Then
			\[
			\qbinom{n}{k}_q=\qbinom{n}{n-k}_q.
			\]
		\end{proposition}
		
		\begin{proof}
			Let $V = \F_q^n$. The map $U \mapsto \{\theta \in V^\star : \theta(U) = 0\}$
			is a bijection between $k$-dimensional subspaces of $V$ and the $(n-k)$-dimensional subspaces of the dual space $V^\star$.
		\end{proof}
		
		One of the difficulties in proving $q$-binomial identities lies in the occurrence of powers of $q$ (for instance, on the right-hand side of Theorem~\ref{thm:q-nan}), which specializes to 1 in the binomial case. The vector space interpretation often explains the appearance of these factor --- in contrast to set complements in the binomial setting, \emph{a complement of a subspace is not unique}.
		We shall see that the complements are counted by a power of $q$.
		
		\begin{lemma}\label{prop:number_of_complements}
			Let $U$ be a $k$-dimensional subspace of $V=\mathbb{F}_q^n$. Then $U$ has $q^{k(n-k)}$ distinct complements inside $V$.
		\end{lemma}
		
		\begin{proof}
			Let $e_1,\ldots, e_n$ be the canonical basis of $V$.
			We may assume without loss of generality that $U$ is spanned by $e_1,\ldots, e_k$
			and identify $U$ with $\F_q^k$. Each complement of $U$ has a unique basis of the form 
			\[
				u_1+e_{k+1}, \ldots, u_{n-k}+e_{n},
			\]
			where $u_1, \ldots, u_{n-k}\in \F_q^k$. 
			The $(n-k)$-tuples $(u_1,\ldots, u_{n-k})$ are in bijection with $(n-k) \times k$ matrices with entries
			in $\mathbb{F}_q$, of which there are
			$q^{k(n-k)}$.
%
%
%
		\end{proof}
		\begin{example}
			As seen in Figure~\ref{fig:f3-example}, when $q = 3$, $n=2$ and $k=1$ each line has three distinct complements in $\F_3^2$.
		\end{example}
		
		\begin{lemma}\label{cor:sub-of-complement}
			Let $U$ be an $r$-dimensional subspace of $V=\F_q^n$. The number of $s$-dimensional subspaces $W$ of $V$ such that $U\cap W=0$ is $q^{r(n-r)}\qbinom{n-r}{s}_q$.
		\end{lemma}
		\begin{proof}
			By Proposition~\ref{prop:q-binom-defn} there are $\qbinom{n-r}{s}_q$ $(r+s)$-dimensional subspaces $X$ of $V$ containing $U$. 
			By Lemma~\ref{prop:number_of_complements} there are $q^{r(n-r)}\qbinom{n-r}{s}_q$ 
			complements $W$ of~$U$ inside each $X$.
		\end{proof}
		
		Another important difference comparing $q$-binomial coefficients and subspaces to
		 binomial coefficients and subsets is the analogue of set union: taking a direct sum of two subspaces \textit{loses} information about the vectors which belong to neither of the direct summands. A more careful approach, frequently used in our proof of Theorem~\ref{thm:q-nan} is through \textit{extensions}, counted bijectively in the following proposition.
		\begin{proposition}\label{prop:number_of_extension}
			Let $U\leq W$ be given $k$- and $\ell$-dimensional subspace of $V=\mathbb{F}_q^n$. The number of $m$-dimensional extensions $E$ of $U$, such that $E\cap W=U$, equals
		\begin{align*}
			q^{(m-k)(\ell-k)}\qbinom{n-\ell}{m-k}_q.
		\end{align*}
	\end{proposition}
		\begin{proof}
			By Lemma~\ref{cor:sub-of-complement}, applied to $V/U$, there are \smash{$q^{(m-k)(\ell-k)}\qbinom{(n-k)-(\ell-k)}{m-k}_q$} $(m-k)$-dimensional subspaces $\overline{E}$ of $V/U$ such that $\overline{E}\,\cap\, W/U=0$. Now again use that such subspaces of the quotient space $V/U$ are in bijection with the $m$-dimensional subspaces $E$ of $V$ such that $U\leq E$ and $E\cap W=U$.
		\end{proof}
%
		\begin{example}
			If $k=\ell$, Proposition~\ref{prop:number_of_extension} correctly simplifies to Lemma~\ref{lem:n-k_choose_l-k}.
			If instead $k=0$ and $m=n-\ell$, Proposition~\ref{prop:number_of_extension} simplifies to  Lemma~\ref{prop:number_of_complements}.
		\end{example}
		
		As a corollary, we deduce a $q$-lift of Vandermonde's convolution. Again it has a bijective proof.
		
		\begin{corollary}[$q$-Vandermonde's convolution]\label{prop:q-Vandermonde2}
			Let $m\leq n$ be non-negative integers. Then
			\begin{align*}
				\qbinom{n}{m}_q=\sum_{k=0}^m q^{k(n-\ell-m+k)}\qbinom{\ell}{k}_q\qbinom{n-\ell}{m-k}_q
			\end{align*}
			for any non-negative integer $\ell\leq n$.
		\end{corollary}
		
		
		
		\begin{proof}
			Let $W$ be a fixed $(n-l)$-dimensional subspace of $V=\F_q^n$. The $m$-dimensional subspaces of $V$ counted by the left-hand side are in one-to-one correspondence with $m$-dimensional extensions $E$ of
			 the $(m-k)$-dimensional subspace $U$ of $W$ defined by $U = E\cap W$. By Proposition~\ref{prop:q-binom-defn} and Proposition~\ref{prop:number_of_extension}, there are $q^{k(n-\ell-m+k)}\qbinom{\ell}{k}_q\qbinom{n-\ell}{m-k}_q$ many such pairs $(E, U)$. Summing over $k$ gives the result. 
		\end{proof}
		
		\begin{remark}\label{remark:primePowersSuffice}
			Since the $q$-integers are polynomials in $q$ with coefficients in the field~$\mathbb{C}$,
			to prove a $q$-binomial identity it suffices to prove it for infinitely many values of $q$,
			such as the prime powers. In particular,
			the identities obtained in this section  hold for all $q \in \mathbb{C}$. 
		\end{remark}
		
		\section{A new bijective proof of Theorem~\ref{thm:q-nan}}\label{sec:new-bijective-proof}
		
		In this section we construct the bijective map that establishes Theorem~\ref{thm:q-nan}  --- the main difficulty 
		is to construct the sets counted by the left- and right-hand side. We observe immediately that the left-hand side of Theorem~\ref{thm:q-nan} counts pairs of subspaces of dimensions $s$ and $t$ in $\mathbb{F}_q^m$ and $\mathbb{F}_q^{m+e}$, respectively. The right-hand side requires a more detailed description. A big step toward this goal is made in Lemma~\ref{lemma:septuples}.
		
		\begin{definition}\label{def:fixed}
			Let $X$ be a subspace of an $\F_q$-vector space $Y$. 

\begin{defnlist}\item			For a non-negative integer $d\leq\dim X$, the \emph{fixed $d$-dimensional subspace of~$X$}, denoted $U_{X, d}$ is defined as follows. Consider the canonical basis of $Y$ and assume the coordinate-wise lexicographical order of its elements (recall that $Y$ is a finite set). We may label them with $1,2,\ldots, q^n$, and choose the $d$ linearly independent elements of $X$ whose set of labels is lexicographically smallest. The vector space $U_{X,d}$ is defined as the span of these $d$ vectors. 
	
	\item		The \emph{fixed complement of $X$ in $Y$}, denoted $X^{\perp Y}$, is the complement which is lexicographically smallest as a finite set of vectors.
	\end{defnlist}
		\end{definition}

		\begin{lemma}\label{lemma:septuples}
			Let $m,s,t,e,j,k$ be non-negative integers. Let $W$ be an $m$-dimensional subspace of $V=\mathbb{F}_q^{m+e}$.  The number of quadruples 
			\begin{align*}
				(A,B,C, D)
			\end{align*}
			of subspaces of $V$ satisfying the following conditions:
			\begin{thmlist}
				\item[$\bullet$] $\dim A=t$, $\dim B=s$, $\dim C=s+t-k$, and $\dim D=j-e$,
				\item[$\bullet$] $A\leq C\leq W$ and $D\leq U_{C, s-e}$,
				\item[$\bullet$] $A\cap U_{C, s-e}=D$ and $B\cap A^{\perp W}=A^{\perp W}\cap C$
			\end{thmlist}
			equals
			\begin{align*}
				\qbinom{m}{s+t-k}_q\cdot\qbinom{s-e}{s-j}_q\cdot q^{(s-j)(t+e-j)}\qbinom{t+e-k}{j-k}_q\cdot q^{k(m-s-t+k)}\qbinom{t+e}{k}_q.
			\end{align*}
		\end{lemma}
		
		The formula and conditions arise naturally in the process of construction presented in the proof below. 
		Since there is a $k$-dimensional subspace of $\F_q^n$ if and only if $0 \le k \le n$, the formula above is non-zero if and only if the values of $m,s,t,e,j,k$ allow all the conditions to be satisfied:
		for instance, $t=\dim A\leq \dim C=s+t-k$, etc.
		
		
		\begin{figure}[h]\centering
			\caption{The construction of quadruples from Lemma~\ref{lemma:septuples}. The rectangles represent chosen subspaces, with the area underneath the dashed line corresponding to the $m$-dimensional subspace $W$ of $V$. }\label{fig:septuples}
			\vspace*{3pt}
%
			\begin{subfigure}[t]{0.5\linewidth}\centering
				\begin{tikzpicture}[scale=1.35]
					\draw[dashed, thick] (0.5, 2.3) -- (3.5, 2.3);
					\fill[mygreen!10] 
					(1, 2.3)--(1, 1)--(3, 1) -- (3,2.3) -- cycle;
					\fill[mygreen!20] (2,1)--(3,1)--(3,2.3)--(2,2.3)--cycle;
					\draw[ultra thick, mygreen] 
					(1, 2.3)--(1, 1)--(3, 1) -- (3,2.3) -- cycle;
					
					\draw[thick] (0.5, 0.5) rectangle (3.5,3);
					

					\node[align=left] at (2, 0.75) {$C$};
					
					\node[align=left] at (3.7, 1.6) {$V$};
					\node[align=left] at (3.3, 0.7) {$W$};
					\node[align=left] at (2.5, 2) {$U_{C, s-e}$};
				\end{tikzpicture}
				\caption{ Choose $C\cong\mathbb{F}_q^{s+t-k}$ inside $W$ \\ \centering
					in $\qbinom{m}{s+t-k}_q$ many ways.}\label{fig:choose-C}
			\end{subfigure}\hfill
%
%
%
%
%
			\begin{subfigure}[t]{0.5\linewidth}\centering
				\begin{tikzpicture}[scale=1.35]
					\fill[mygreen!10] 
					(1, 2.3)--(1, 1)--(3, 1) -- (3,2.3) -- cycle;
					\fill[mygreen!20] (2,1)--(3,1)--(3,2.3)--(2,2.3)--cycle;
					\fill[mygreen!30] (2,1)--(3,1)--(3, 1.6)--(2, 1.6)--cycle;
					\draw[thick] (0.5, 0.5) rectangle (3.5,3);
					\draw[dashed, thick] (0.5, 2.3) -- (3.5, 2.3);
					\draw[ultra thick, mygreen] (2, 1.6)--(2, 1)--(3, 1) -- (3,1.6)--cycle;

					\node[align=left] at (2, 0.75) {$C$};
					\node[align=left] at (2.5, 1.3) {$D$};
					\node[align=left] at (3.7, 1.6) {$V$};
					\node[align=left] at (3.3, 0.7) {$W$};
					\node[align=left] at (2.5, 2) {$U_{C, s-e}$};
					
				\end{tikzpicture}
				\caption{ Choose $D\cong\mathbb{F}_q^{j-e}$ inside $U_{C, s-e}$ \\ \centering in $\qbinom{s-e}{s-j}_q$ many ways.}\label{fig:choose-D}
			\end{subfigure}\\
			
			\begin{subfigure}[t]{0.5\linewidth}\centering
				\begin{tikzpicture}[scale=1.35]
					\fill[mygreen!10] 
					(1, 2.3)--(1, 1)--(3, 1) -- (3,2.3) -- cycle; 
					\fill[mygreen!20] (2,1)--(3,1)--(3,2.3)--(2,2.3)--cycle; 
					\fill[pattern = north west lines, pattern color = mygreen!40] (0.5,0.5)--(3.5,0.5)--(3.5,2.3)--(0.5,2.3)--cycle;
					\fill[mygreen!30] (2,1)--(3,1)--(3, 1.6)--(2, 1.6)--cycle; 
					\fill[mygreen!40] (1,1)--(2,1)--(2,1.6)--(1,1.6)--cycle; 
					\draw[thick] (0.5, 0.5) rectangle (3.5,3);
					\draw[dashed, thick] (0.5, 2.3) -- (3.5, 2.3);
					\draw[ultra thick, mygreen] (1, 1.6)--(1, 1)--(3, 1) -- (3,1.6)--cycle;

					\node[align=left] at (1.9, 1.3) {$A$};
					\node[align=left] at (2, 0.75) {$C$};
					\node[align=left] at (2.5, 1.3) {$D$};
					\node[align=left] at (3.7, 1.6) {$V$};
					\node[align=left] at (3.3, 0.7) {$W$};
					\node[align=left] at (2.5, 2) {$U_{C, s-e}$};
					\node[align=left] at (0.95, 1.95) {$A^{\perp W}$};
				\end{tikzpicture}
				\caption{ \centering Extend $D$ to $A\cong\mathbb{F}_q^{t}$ 
					in~$C$\\so that $A\cap U_{C, s-e}=D$~ in\\$q^{(s-j)(t+e-j)}\qbinom{t+e-k}{j-k}_q$~many~ways.}\label{fig:extend-D}
			\end{subfigure}\hfill
%
%
%
%
			\begin{subfigure}[t]{0.5\linewidth}\centering
				\begin{tikzpicture}[scale=1.35]
					\fill[mygreen!10] 
					(1, 2.3)--(1, 1)--(3, 1) -- (3,2.3) -- cycle; 
					\fill[mygreen!20] (2,1)--(3,1)--(3,2.3)--(2,2.3)--cycle; 
					\fill[pattern = north west lines, pattern color = mygreen!40] (0.5,0.5)--(3.5,0.5)--(3.5,2.3)--(0.5,2.3)--cycle;
					\fill[mygreen!30] (2,1)--(3,1)--(3, 1.6)--(2, 1.6)--cycle; 
					\fill[mygreen!40] (1,1)--(2,1)--(2,1.6)--(1,1.6)--cycle; 
					\draw[thick] (0.5, 0.5) rectangle (3.5,3);
					\draw[dashed, thick] (0.5, 2.3) -- (3.5, 2.3);
					\draw[ultra thick, mygreen] (1, 1.6)--(1, 1)--(3, 1) -- (3,1.6)--cycle;
					\fill[pattern = north east lines, pattern color = mygreen!70] (1,1.45)--(3,1.45)--(3, 2.8)--(1,2.8)--cycle;
					\draw[ultra thick, mygreen] (1,1.45)--(3,1.45)--(3, 2.8)--(1,2.8)--cycle;

					\node[align=left] at (2, 2.5) {$B$};
					\node[align=left] at (1.9, 1.3) {$A$};
					\node[align=left] at (2, 0.75) {$C$};
					\node[align=left] at (2.5, 1.3) {$D$};
					\node[align=left] at (3.7, 1.6) {$V$};
					\node[align=left] at (3.3, 0.7) {$W$};
					\node[align=left] at (2.5, 2) {$U_{C, s-e}$};
					\node[align=left] at (0.95, 1.95) {$A^{\perp W}$};
				\end{tikzpicture}
				\caption{ \centering Extend $A^{\perp W}\cap C$ to $B\cong\mathbb{F}_q^s$  in~$V$\\so that $B\cap A^{\perp W}=A^{\perp W}\cap C$ in\\$q^{k(m-s-t+k)}\qbinom{t+e}{k}_q$~many~ways.}\label{fig:extend-A^C-2}
			\end{subfigure}

		\end{figure}

		\begin{proof}[Proof of Lemma~\ref{lemma:septuples}]
			Choose an $(s+t-k)$-dimensional subspace $C$ of $W$ (Figure~\ref{fig:choose-C}) in $\qbinom{m}{s+t-k}_q$ many ways and consider its fixed $(s-e)$-dimensional subspace $U_{C, s-e}$ defined in Definition~\ref{def:fixed}(a). Next, choose a $(j-e)$-dimensional subspace $D$ of $U_{C, s-e}$ (Figure~\ref{fig:choose-D}) in $\qbinom{s-e}{j-e}_q=\qbinom{s-e}{s-j}_q$ many ways, where the equality uses Proposition~\ref{lemma:n-k}. By Proposition~\ref{prop:number_of_extension}, $D$ can be extended to a $t$-dimensional subspace $A$ of $C$ (Figure~\ref{fig:extend-D}), such that $A\cap U_{C, s-e}=D$, in 
			\[ q^{(s-j)(t+e-j)}\qbinom{t+e-k}{t+e-j}_q=q^{(s-j)(t+e-j)}\qbinom{t+e-k}{j-k}_q \] many ways; 
			again the equality uses Proposition~\ref{lemma:n-k}.
			Let $A^{\perp W}$ be the fixed complement of $A$ in $W$ defined in Definition~\ref{def:fixed}(b). Then by Proposition~\ref{prop:number_of_extension}, $A^{\perp W}\cap C\leq A^{\perp W}$ 
			can be extended to an
			$s$-dimensional vector space $B$ in $V$ (Figure~\ref{fig:extend-A^C-2}), such that $B\cap A^{\perp W}=A^{\perp W}\cap C$, in $q^{k(m-s-t+k)}\qbinom{t+e}{k}_q$ many ways. 
			
			We have therefore chosen a quadruple satisfying the specified conditions, and the total number of such quadruples is the product of the number of choices made at each step; this gives the desired expression.
		\end{proof}
		To prove Theorem~\ref{thm:q-nan} when $e < 0$, we need the following standard result.
		
		\begin{lemma}\label{lemma:specialization}
			Let $f \in \mathbb{Q}(q)[X]$ have degree $d$ as a polynomial in $X$. If $f(q^h)=0$ for at least $d+1$ values of $h\in\mathbb{N}_0$, then $f=0$. Moreover, if $g \in \Q(q)[X,X^{-1}]$ vanishes at $(q^h, q^{-h})$ for infinitely many values of $h\in\N_0$, then $g=0$.
		\end{lemma}
		
		\begin{proof}
			The first part is a well-known fact about polynomials with coefficients in an integral domain. For the second part, take $a$
			sufficiently large that~$X^a g(X, X^{-1})$ is a polynomial in $X$. The conclusion follows from the first part.
		\end{proof}
		
		\begin{proof}[Proof of Theorem~\ref{thm:q-nan} \emph{($q$-Pfaff--Saalsch\"utz identity)}]
			Assume $0\leq e\leq s$ and let $q$ be a prime power. Let $W=U_{V, m}$ be the fixed $m$-dimensional subspace of $V=\mathbb{F}_q^{m+e}$. 
			
			Consider the 
			set $\mathcal{T}$ of pairs $(A,B)$ of subspaces of $V$ satisfying $\dim  A=t$, $\dim B=s$, and $A\leq W$. By Proposition~\ref{prop:q-binom-defn}, $|\mathcal{T}|$ is the left-hand side of the identity in Theorem~\ref{thm:q-nan}.
			
			Let $\mathcal{S}_{j,k}$ denote the set of quadruples described in Lemma~\ref{lemma:septuples} for given $j,k\in\N_0$. Define $\mathcal S=\bigcup_{j,k\in\mathbb{N}_0}\mathcal{S}_{j,k}$. We apply the $q$-Vandermonde's convolution (Corollary~\ref{prop:q-Vandermonde2}) to the term $\qbinom{m+j}{s+t}_q$ of the right-hand side, 
			thus introducing a new summation,
			followed by Proposition~\ref{lemma:subset_of_subset} to the product $\qbinom{t+e}{j}_q\qbinom{j}{k}_q$:
			\begin{align*}
				&\sum_{j\geq 0}q^{(s-j)(t+e-j)}\qbinom{t+e}{j}_q\qbinom{s-e}{s-j}_q\qbinom{m+j}{s+t}_q\\
				&=\sum_{j\geq 0}q^{(s-j)(t+e-j)}\qbinom{t+e}{j}_q\qbinom{s-e}{s-j}_q\sum_{k\geq0}q^{k(m-s-t+k)}\qbinom{m}{s+t-k}_q\qbinom{j}{k}_q\\
				&=\sum_{j\geq0}\sum_{k\geq0}q^{(s-j)(t+e-j)}q^{k(m-s-t+k)}\qbinom{t+e}{k}_q\qbinom{t+e-k}{j-k}_q\qbinom{s-e}{s-j}_q\qbinom{m}{s+t-k}_q.
			\end{align*}
			Therefore, since the sets $\mathcal{S}_{j,k}$ are disjoint, by Lemma~\ref{lemma:septuples}, $|\mathcal{S}|$ equals the right-hand side of the identity from Theorem~\ref{thm:q-nan}. 
			
%
			
			It remains to exhibit a bijection between the sets $\mathcal{T}$ and $\mathcal{S}$. One direction of the bijection is the natural projection $\mathcal{S} \rightarrow \mathcal{T}$ defined by
			\[  (A,B,C, D) \mapsto(A,B).\]
			To construct the inverse we retrace the proof of Lemma~\ref{lemma:septuples}. Define $C=A\,\oplus\, (B\,\cap\, A^{\perp W})$ and $D= A\,\cap\, U_{A\,\oplus\, (B\,\cap\, A^{\perp W}), s-e}$. By construction, the conditions listed in the said lemma are satisfied for such quadruple.
			Therefore, the map $\mathcal{T} \rightarrow \mathcal{S}$ defined by
			\[  
				(A,B)\mapsto\big(A,B,A\oplus (B\cap A^{\perp W}), 
				A\cap U_{A\oplus (B\cap A^{\perp W}), s-e}\big)
			 \]
			is the desired inverse. 
%
			This proves Theorem~\ref{thm:q-nan}  for $0 \le e \le s$
			when $q$ is a prime power, and the result for arbitrary $q$ follows
			from Remark~\ref{remark:primePowersSuffice}.
			
			To prove the theorem for $-t \le e \le 0$, we
			define a function $f$ by 
		\[ \begin{split}
				f(X)&=\qbinom{m}{t}_q\frac{(q^{m-s+1}X;q)_s}{(q;q)_s}\\
				& \quad -\sum_{j\geq0}q^{(s-j)(t-j)}\frac{(q^{t-j+1}X;q)_j}{(q;q)_j}\cdot\frac{X^{s-j}(q^{j+1}X^{-1};q)_{s-j}}{(q;q)_{s-j}}\cdot\qbinom{m+j}{s+t}_q\end{split}
\]
			where $(a;q)_\alpha$ denotes the shifted factorial product
			\begin{align*}
				(a;q)_\alpha=(1-a)(1-qa)\dots(1-q^{\alpha-1}a).
			\end{align*}
			Note that $\qbinom{n}{k}_q=(q^{n-k+1};q)_k/(q;q)_k$ and hence \[ f(q^e) = \qbinom{m}{t}_q\qbinom{m+e}{s}_q - \sum_{j \ge 0} q^{(s-j)(t-j)} \qbinom{t+e}{j}_q q^{e(s-j)} \qbinom{s-e}{j}_q \qbinom{m+j}{s+j}_q \]
			which is the difference of the two sides in Theorem~\ref{thm:q-nan}.
			By the combinatorial argument above we have $f(q^e) = 0$ for all $e\in\{0,1,\ldots,s\}$.
			Since $f(X)\in\mathbb{Q}(q)[X]$ and $\deg f\leq s$, Lemma~\ref{lemma:specialization} implies that $f=0$. 
			In particular, $f(q^e)=0$ for all $e \in \Z$ such that $-t \le e \le s$.
		\end{proof}
		
		\begin{corollary}[Pfaff--Saalsch\"utz identity]\label{cor:nan}
			Let $m,s,t\in\N_0$. If $e\in\Z$ and \break$-t \le e \le s$ then
			\begin{align*}
				\binom{m}{t}\binom{m+e}{s}=\sum_{j\geq 0}\binom{t+e}{j}\binom{s-e}{s-j}\binom{m+j}{s+t}.
			\end{align*}
		\end{corollary}
		
		\begin{proof}
			Take $q\rightarrow 1$ in Theorem~\ref{thm:q-nan}.
		\end{proof}
		
		\begin{remark}
			In the binomial case, one may consider $A,B,C,D,V,W$ as sets and subsets instead of vector spaces and subspaces. Replacing $A^{\perp W}$ with the set difference $W\setminus A$, the steps in the proof of Theorem~\ref{thm:q-nan} simplify significantly, giving a new bijective proof of the binomial identity in Corollary~\ref{cor:nan}.
		\end{remark}
		
		To illustrate the generality of Theorem~\ref{thm:q-nan}, we show it is equivalent to Stanley's identity
		\cite{Sta70};
		this makes clear a symmetry that is hidden in our statement of Theorem~\ref{thm:q-nan}.
		
		\begin{example}[Stanley's $q$-binomial identity]\label{ex:Stanley}
			Let $x,y,A,B\in\N$. Then
			\begin{align*}
				\qbinom{x+A}{B}_q\qbinom{y+B}{A}_q=\sum_{K\geq 0}q^{(A-K)(B-K)}&\qbinom{x+y+K}{K}_q\qbinom{y}{A-K}_q\qbinom{x}{B-K}_q
			\end{align*}
			is equivalent to Theorem~\ref{thm:q-nan} via the simultaneous substitutions
			\begin{align*}
				(m,e,s,t,j)&\mapsto(B+y,A+x-B-y,A+x-B, B+y-A, K-B+x),\\
				(A,B,x,y,K)&\mapsto(m-t, m+e-s, t+e, s-e, m-s-t+j).
			\end{align*}
		\end{example}
		
		Another notable symmetric form of the $q$-Pfaff--Saalsch\"utz identity was proved by Zeilberger
		\cite{Zei87} lifting Foata's proof \cite{Foa65} of the corresponding binomial identity.
		
		\begin{example}[Foata/Zeilberger form of the $q$-Pfaff--Saalsch\"utz identity]\label{ex:Zeil}
			Let $a, b, c,\break k, n\in \N_0$ with $k$ the smallest. Then setting
			\[
				(m,e,s,t,j)\mapsto (b+c, a-c, a+k, c-k, a-n)
			\]
			in Theorem~\ref{thm:q-nan} gives
			\[
				\qbinom{b+c}{c-k}_q\qbinom{b+a}{a+k}_q=\sum_{n\ge 0}q^{(n-k)(n+k)}\qbinom{a-k}{a-n}_q\qbinom{c+k}{n+k}_q\qbinom{a+b+c-n}{a+c}_q.
			\]
			Now multiply both sides by $\qbinom{c+a}{c+k}_q$, use the definition of $q$-binomial coefficient~\eqref{eq:defn}, and apply Proposition~\ref{lemma:n-k} to obtain (1) in \cite{Zei87}:
			\begin{align*}
				\qbinom{a+b}{a+k}_q\qbinom{b+c}{b+k}_q&\qbinom{c+a}{c+k}_q\\&=\sum_n q^{n^2-k^2}\frac{[a+b+c-n]_q!}{[a-n]_q![b-n]_q![c-n]_q![n+k]_q![n-k]_q!}.
			\end{align*}
			where the sum is over all $n$ such that $k \le n \le \min(a,b,c)$.
		\end{example}
		
		As we mentioned in the introduction, multiple $q$-identities of this form have been 
		rediscovered over the years. In the language of hypergeometric series it becomes clear that they are all equivalent to the $q$-Pfaff--Saalsch\"utz formula. For an introduction that uses the notation
		below, see \cite[Ch.~3]{Sla66}.
		
		\begin{remark}\label{remark:hypergeometric}
			The $q$-binomial identity in  Theorem~\ref{thm:q-nan} is equivalent to the $q$-Pfaff--Saalsch\"utz formula
			\begin{align*}
				\setlength\arraycolsep{1.3pt}
				{}_3 {\Phi}_2\!\left(\begin{matrix}q^{s-m-e}, q^{-t-e}, q^{-t}\\q^{-m-t-e}, 
					q^{s+1-t-e}\end{matrix}~;q,q\right)=\frac{(q^{-m};q)_t(q^{-s-t};q)_t}{(q^{-m-t-e};q)_t(q^{-s+e};q)_t},
			\end{align*}
			which in full generality \cite[(5)]{Zen89} states that for all $a,b,c\in\mathbb{C}$ and $n\in\N_0$:
			\begin{align*}
				\setlength\arraycolsep{2pt}
				{}_3 {\phi}_2\!\left(\begin{matrix}a, b, q^{-n}\\c, 
					q^{1-n}ab/c\end{matrix}~;q,q\right)=\frac{(c/a;q)_n(c/b;q)_n}{(c;q)_n(c/ab;q)_n}.
			\end{align*}
			Conversion between $q$-binomial and hypergeometric identities is essentially routine:
			the procedure is exemplified in the binomial case in \hbox{\cite[\S 5.5]{GKP89}.}
		\end{remark}
		
		\section{The quantum Pfaff--Saalsch\"utz identity}\label{sec:q-to-quantum}
		Let $q\in\C\setminus\{0, \pm1\}$. To lift the identity in Theorem~\ref{thm:q-nan} to Theorem~\ref{thm:multiplication_rule}, we 
		first express its quantum analogue. 
		We define the \textit{quantum integer} by
		\begin{align*}
			\{n\}_q=q^{n-1}+q^{n-3}+\dots+q^{1-n}=\frac{q^n-q^{-n}}{q-q^{-1}}
		\end{align*}
		for $n \in \N_0$. 
		The definitions of \emph{quantum factorial} and \emph{quantum binomial coefficient} again follow naturally:
		$\{n\}_q!=\{n\}_q\{n-1\}_q\ldots\{1\}_q$ and
		\begin{align*}
			\Qbinom{n}{k}_q=\frac{\{n\}_q!}{\{k\}_q!\{n-k\}_q!},
		\end{align*}
		for $n, k \in \N_0$. We set 
		$\Qbinom{n}{k}_q=0$ whenever $n < k$. 
		As we have seen, the \break$q$-binomial coefficients admit multiple combinatorial interpretations. On the other hand, quantum binomial coefficients are more often seen in algebra: for example
		$\Qbinom{n}{k}$ is the character of the representation $\bigwedge^k \Sym^{n-1} \C^2$
		of $\SL_2(\C)$. The connection between the 
		two, established in the following two results, 
		is a fundamental bridge in algebraic combinatorics, connecting combinatorial enumeration, plethysms of symmetric
		functions, and representation theory.
		
		\begin{lemma}\label{lemma:qint-to-Qint}
			Let $n\in\N_0$. Then:
			\begin{align*}
				[n]_{q^2}=q^{n-1}\{n\}_q.
			\end{align*}
		\end{lemma}
		
		\begin{proof}
			This follows directly from the definitions:
			\[
			[n]_{q^2}=\frac{q^{2n}-1}{q^2-1}=\frac{q^n}{q}\cdot\frac{q^n-q^{-n}}{q-q^{-1}}=q^{n-1}\{n\}_q.\qedhere
			\]
		\end{proof}

		\begin{corollary}\label{cor:q-to-quantum}
			Let $k\leq n$ be non-negative integers. Then
			\[
			\qbinom{n}{k}_{q^2}=q^{k(n-k)}\Qbinom{n}{k}_q.
			\]
		\end{corollary}
		
		\begin{proof}
			This is immediate from Lemma~\ref{lemma:qint-to-Qint}.
		\end{proof}
		
		We can now state and prove the quantum version of Theorem~\ref{thm:q-nan}.
		
		\begin{theorem}[Quantum Pfaff--Saalsch\"utz identity]\label{thm:quantum-nan}
			Let $m,s,t\in\N_0$. If $e\in \Z$ and $-t \le e \le s$ then
			\begin{align*}
				\Qbinom{m}{t}_q\Qbinom{m+e}{s}_q=\sum_{j\geq 0}\Qbinom{t+e}{j}_q\Qbinom{s-e}{s-j}_q\Qbinom{m+j}{s+t}_q.
			\end{align*}
		\end{theorem}
		
		\begin{proof}[Proof of Theorem~\ref{thm:quantum-nan}]
			
			By Corollary~\ref{cor:q-to-quantum}, replacing $q$ with $q^2$ in Theorem~\ref{thm:q-nan} gives
			the required identity up to a power of $q$. The routine calculation 
			\begin{align*}
				q^{t(m-t)}q^{s(m+e-s)}=q^{2(s-j)(t+e-j)}q^{j(t+e-j)}q^{(s-j)(j-e)}q^{(s+t)(m+j-s-t)}
			\end{align*}
			shows that the powers of $q$ cancel on both sides.
		\end{proof}
		
		The form of Theorem~\ref{thm:quantum-nan} we need for the proof of Theorem~\ref{thm:multiplication_rule} is as follows.
		\begin{corollary}
			Let $h,s,t\in\N_0$. If $b,c\in \Z$, $t-c+b\geq 0$, and $s-b+c\geq0$ then
			\begin{align*}\label{eq:quantumPS}
				\Qbinom{h+c}{t}_q\Qbinom{h+b}{s}_q=\sum_{i\geq 0}\Qbinom{t-c+b}{i-c}_q\Qbinom{s-b+c}{i-b}_q\Qbinom{h+i}{t+s}_q.
			\end{align*}
		\end{corollary}
		\begin{proof}
			Make the substitution $(m,e,j)\mapsto(h+c, b-c, i-c)$ in Theorem~\ref{thm:quantum-nan}. 
			The top entries in Theorem~\ref{thm:quantum-nan} are non-negative by assumption, 
			and they remain non-negative after the substitution.
		\end{proof}

		\section{Lusztig's integral form of \texorpdfstring{$\U_q(\sl_2)$}{Uq(sl2)}}\label{sec:quantum}
		The quantum group $\U_q(\sl_2)$ is the unital associative algebra over $\Q(q)$ with generators $E$, $F$, $K$, $K^{-1}$ subject to the 
		relation
		\begin{equation*}
			\begin{aligned}
				& KEK^{-1} &=q^2E,\quad & KFK^{-1}&=q^{-2}F,\quad & EF-FE&= \frac{K-K^{-1}}{q-q^{-1}}
			\end{aligned}
		\end{equation*}
		and the expected $KK^{-1} = K^{-1} K = 1$. By
		Theorem~3.1.5 in \cite{HK02} it has a
		triangular decomposition as
		\[\U_q(\sl_2) = \U_q^-(\sl_2)\otimes_{\Q(q)}\U_q^0(\sl_2)\otimes_{\Q(q)}\U_q^+(\sl_2)\]
		where $U^-_q(\sl_2)$ (resp.~$\U_q^0(\sl_2)$, resp. $\U_q^+(\sl_2)$) is the $\Q(q)$-subalgebra generated by $F$ (resp.~$K^{\pm 1}$, resp.~$E$).
		Let 
		\begin{equation}\label{eq:Ka} [K;a]= \frac{q^aK-q^{-a}K^{-1}}{q-q^{-1}}.\end{equation}
		In \cite{Lus88}, Lusztig defines the \emph{divided powers} $E^{(n)}=\frac{E^n}{\{n\}_q!}$ and $F^{(n)}=\frac{F^n}{\{n\}_q!}$, as well as 
		the \emph{Lusztig elements}
		$\qbinom{K;\hskip1pt c}{t}$ for $t \in \N_0$ and $c\in \Z$ already seen in~\eqref{formula:q-kostant}. 
		By definition $\qbinom{K;\hskip1pt c}{0}=1$
		and to simplify notation we set
		$\qbinom{K}{t}=\qbinom{K;\hskip1pt 0}{t}$.  Lusztig's element $\qbinom{K;c}{t}$ should be thought of as a quantization of the element
		\[\binom{H+c}{t}=\frac{(H+c)(H+c-1)\cdots (H+c-t+1)}{t!},\] 
		belonging to Kostant's $\Z$-form of the enveloping algebra of $\sl_2$.
		Recall in the following definition that $\A=\Z[q,q^{-1}]$. 
		
		\begin{definition}\label{defn:Lusztig}
			\emph{Lusztig's integral form} $\U_\A(\sl_2)$ for $\U_q(\sl_2)$ is the $\A$-algebra generated by the divided powers 
			\smash{$E^{(n)}=\frac{E^n}{\{n\}_q!}$} and \smash{$\frac{F^n}{\{n\}_q!}$} for $n \in \N$ and the elements $K,K^{-1}$ and $\qbinom{K}{t}$ for $t\in \N$.
		\end{definition}
		
		The divided powers satisfy the following relations for $n, m \in \N$:
		\begin{align}
			E^{(n)}E^{(m)} &= \Qbinom{n+m}{n}_q E^{(n+m)} \label{eq:Erelation}, \\ 
			F^{(n)}F^{(m)} &= \Qbinom{n+m}{n}_q F^{(n+m)} \label{eq:Frelation}.
			\intertext{We also have the following relations \cite[\S4.3]{Lus88}:}
			\qbinom{K;c}{t} E^{(n)} &= E^{(n)} \qbinom{K;c+2n}{t} \label{eq:EKrelation}, \\ 
			\qbinom{K;c}{t} F^{(n)} &= F^{(n)} \qbinom{K;c-2n}{t} \label{eq:FKrelation}, \\ 
			E^{(n)} F^{(m)} &= \sum_{t \geq 0} F^{(m-t)} \qbinom{K;2t-m-n}{t} E^{(n-t)} \label{eq:EFrelation}.
		\end{align}
		
		\begin{proposition}[Lusztig]\label{thm:lusztigbasis}
			The $\A$-algebras $\U_\A(\sl_2)$ and $\mathcal{U}^0_\A(\sl_2)$ are free as $\A$-modules.
			\begin{enumerate}
				\item[\emph{(a)}]
				The $\A$-algebra $\U_\A(\sl_2)$ has an $\A$-basis given by the elements \[F^{(a)}K^\delta\qbinom{K}{t}E^{(b)}\]
				for $a,b,t \in \N_0$ and $\delta\in \{0,\min(1,t)\}$.
				\item[\emph{(b)}] The elements $K^\delta\qbinom{K}{t}$
				for $t \in \N_0$ and $\delta\in \{0,\min(1,t)\}$ form a basis for the Cartan subalgebra $U^0_\A(\sl_2)$.\label{thm:lusztigbasis:2}
			\end{enumerate}
		\end{proposition} 
		
		\begin{proof}
			The first part follows from Theorem 4.5 and Proposition 2.17 in \cite{Lus90}. The second part follows from Theorem 4.5 and Proposition 2.14 in \cite{Lus90}.
		\end{proof}
		
		We  now introduce a new $\A$-basis of Lusztig's integral form and its Cartan subalgebra. In the following section, 
		we  combine these results with the multiplication rule from Theorem~\ref{thm:multiplication_rule} to deduce a new presentation of these algebras.
		
		\begin{proposition}\label{prop:binomialbasis}
			Lusztig's integral form and its Cartan subalgebra admit the following bases:
			\begin{enumerate} 
				\item[\emph{(1)}]
				The $\A$-algebra $U^0_\A(\sl_2)$ has an $\A$-basis \[\B = \left\{\qbinom{K}{t} : t\ge 0\right\} \cup \left\{\qbinom{K;1}{t} : t\ge 1\right\};\]
				\item[\emph{(2)}] The elements $F^{(a)}\qbinom{K;\hskip1pt c}{t}E^{(b)}$ for $a,b,t \in \N_0$ and $c\in \{0,\min(1,t)\}$, form an 
				$\A$-basis of $U_{\A}(\sl_2)$.
			\end{enumerate}
		\end{proposition}
		\begin{proof}
			First we compute using~\eqref{eq:Ka} that
			\begin{align*}
				q^t[K;1] - q^{-1}[K;-t] &= q^t \frac{qK - q^{-1}K^{-1}}{q - q^{-1}} - q^{-1} \frac{q^{-t}K - q^tK^{-1}}{q - q^{-1}} \\
				&= \frac{q^{t+1}K - q^{-(t+1)}K}{q - q^{-1}}\\
				&= \{t+1\}_qK.
			\end{align*}
			Multiplying both sides by $[K;0]\cdots[K;-t+1]$ and dividing by $\{t+1\}_q!$, the right-hand side becomes $K\qbinom{K;\hskip1pt 0}{t}$, while the left-hand side becomes $q^t\qbinom{K;1}{t+1}-q^{-1}\qbinom{K;\hskip1pt 0}{t+1}$. Thus
			\begin{align*}
				K\qbinom{K}{t}&=q^{t}\qbinom{K;1}{t+1}-q^{-1}\qbinom{K}{t+1}.
			\end{align*}
			Since integer powers of $q$ are units in $\A$, this gives an invertible change of basis from the basis claimed in (1) to the $\A$-basis of $\U_\A^0(\sl_2)$ in Proposition~\ref{thm:lusztigbasis}(b). This proves (1) and (2) is a direct consequence of (1) together with~Proposition~\ref{thm:lusztigbasis}(a).
		\end{proof}
		\begin{remark}
			It is a notable feature of the basis $\B$ in Proposition~\ref{prop:binomialbasis}(1) that it contains neither~$K$ nor $K^{-1}$.
			Instead, as a special case of Proposition~\ref{prop:binomialbasis}, we may write
			\begin{align*}
				K&=\qbinom{K;1}{1}-q^{-1}\qbinom{K}{1},\\
				K^{-1}&=\qbinom{K;1}{1}-q\qbinom{K}{1}.
			\end{align*}
			More generally Proposition~\ref{prop:binomialbasis} implies 
			Lusztig's $\A$-basis of $U^0_\A(\sl_2)$ 
			defined in Definition~\ref{defn:Lusztig} and
			our new basis $\B$ are related by an invertible change-of-basis
			matrix $M$ such that both $M$ and $M^{-1}$ have coefficients in $\A = \Z[q,q^{-1}]$.
		\end{remark}

		\section{The multiplication rule of Lusztig's elements \texorpdfstring{$\qbinom{K;\hskip1pt c}{t}$}{[K;c // t]} and a~new presentation of Lusztig's integral form \texorpdfstring{$\U_\A(\sl_2)$}{UA(sl2)}}\label{sec:quantum-new}
		
		In this section we present further new results about Lusztig's integral form and its Cartan subalgebra, including the multiplication rule from Theorem~\ref{thm:multiplication_rule}, and a new presentation of these algebras.
		
		\begin{proof}[Proof of Theorem~\ref{thm:multiplication_rule}]
			We first unpack 
			the definition of Lusztig's elements in~\eqref{formula:q-kostant}:
			\begin{align*}
				\qbinom{K;c}{t}=\frac{[K;c][K;c-1]\cdots [K;c-t+1]}{\{t\}_q!}=\frac{\prod_{s=0}^{t-1}(Kq^{c-s}-K^{-1}q^{-(c-s)})}{\prod_{s=1}^t(q^s-q^{-s})}.
			\end{align*}
			Observe that upon specializing $K, K^{-1}$ to $q^h, q^{-h}$, respectively, Lusztig element specializes to
			a quantum binomial coefficient:
			\begin{align*}
				\qbinom{K; c}{t}\Bigg|_{\substack{{K=q^h\quad\phantom{x}}\\{K^{-1}=q^{-h}}}}\!\!\!\!=\frac{\prod_{s=0}^{t-1}(q^hq^{c-s}-q^{-h}q^{-(c-s)})}{\prod_{s=1}^t(q^s-q^{-s})}=\Qbinom{h+c}{t}_q.
			\end{align*}
			In particular, for any integer $h\geq \max\{-b, -c\}$, the multiplication rule claimed
			in this theorem specializes to the identity in Theorem~\ref{thm:quantum-nan}. Since there are infinitely many such values of $h$, we conclude the proof by applying Lemma~\ref{lemma:specialization} to the function $g\in\Q(q)[K,K^{-1}]$ defined by
			\[
			g(K,K^{-1})=\qbinom{K; c}{t}\qbinom{K; b}{s}-\sum_{i\ge 0} \Qbinom{t-c+b}{i-c}_q\Qbinom{s-b+c}{i-b}_q\qbinom{K;i}{t+s},
			\]
			treating $K$ as an indeterminant.
		\end{proof}
		
		To give the new presentations, we now need one more identity in Lusztig's integral form, which we state and prove in the proposition below.
		Recall that $\A = \Z[q,q^{-1}]$ and note that the $\A$-algebra $U^0_\A(\sl_2)$ has an automorphism $\Phi_{c}$ for each $c\in \Z$, defined 
		by restricting the automorphism of $\Q(q)[K,K^{-1}]$ which sends $K\mapsto q^{c}K$ and $K^{-1}\mapsto q^{-c}K^{-1}$.
		In the following proposition, we set $\qbinom{K;\hskip1pt c}{t}=0$ if $t<0$.

		\begin{proposition}\label{prop:cplus2}
			Let $t\ge 1$ and $c\in \Z$. The following relation holds 
			\begin{equation*}
				\qbinom{K;c+2}{t}=(q^t+q^{-t})\qbinom{K;c+1}{t}-\qbinom{K;c}{t}+\qbinom{K;c}{t-2}.
			\end{equation*}
		\end{proposition}
		\begin{proof}
			First observe that this relation is obtained by applying the automorphism $\Phi_c$ to the following relation:
			\begin{equation}\label{eq:zeroplus2}
				\qbinom{K;2}{t}=(q^t+q^{-t})\qbinom{K;1}{t}-\qbinom{K}{t}+\qbinom{K}{t-2}.
			\end{equation}
			We prove this simplified identity by considering two cases. 
			We give full details as it is easy to take an inefficient path through these calculations.
			\begin{itemize}
				\item Case $t=1$: by~\eqref{eq:Ka} we have
				\begin{align*}
					\quad (q+q^{-1})\qbinom{K;1}{1}-\qbinom{K;0}{1} &=(q+q^{-1})\frac{qK-q^{-1}K^{-1}}{q-q^{-1}}-\frac{K-K^{-1}}{q-q^{-1}}\\
					&= \frac{q^2K-q^{-2}K^{-1}}{q-q^{-1}}\\
					&= \qbinom{K;2}{1}.
				\end{align*}

				\item Case $t\ge 2$:
				Multiplying~\eqref{eq:zeroplus2} by $\{t\}_q! \bigl/ [K;0]\cdots [K;2-t+1]$ 
				the left-hand side becomes $[K;2][K;1]$ and the right-hand side becomes
				$(q^t+q^{-t})[K;1][K;2-t]  - [K;2-t][K;1-t]+\{t\}_q\{t-1\}_q$, 
				so it suffices to show that
				\[\quad [K;2][K;1] - (q^t+q^{-t})[K;1][K;2-t] + [K;2-t][K;1-t]\]
				equals $\{t\}_q\{t-1\}_q$.
				Multiplying further by $(q-q^{-1})^2$, we obtain
				\begin{align} (q^2K&-q^{-2}K)(qK-q^{-1}K^{-1}) \nonumber \\ &\quad - (q^t+q^{-t})(qK-q^{-1}K^{-1})(q^{2-t}K-q^{-2+t}K^{-1}) \nonumber
					\\ &\quad + (q^{2-t}K-q^{t-2}K^{-1})(q^{1-t}K-q^{-1+t}K^{-1})\label{eq:Kproduct} \end{align}
				which we must show equals 
				\[\qquad  (q-q^{-1})^2\{t\}_q\{t-1\}_q=(q^{t}-q^{-t})(q^{t-1}-q^{1-t})
				= q^{2t-1}+q^{1-2t}-q-q^{-1}
				.\]
				The powers of $K$ that may appear in the expansion of~\eqref{eq:Kproduct} are $K^2$, $K^0$ and $K^{-2}$.
				The coefficient of $K^2$ is
				$q^3-(q^t+q^{-t})q^{3-t}+q^{3-2t}=0$.
				Similarly, the coefficient of $K^{-2}$ is
				$q^{-3}-(q^t+q^{-t})q^{-3+t}+q^{-3+2t}=0$.
				Finally, the coefficient of $K^0$ is
				\begin{align*}\qquad 
					&-q^{2-1}-q^{-2+1}+(q^t+q^{-t})(q^{1+(-2+t)}+q^{-1+(2-t)})\\ &\qquad -q^{(2-t)+(-1+t)} -q^{(t-2)+(1-t)} = q^{2t-1}+q^{1-2t}-q-q^{-1}
				\end{align*}
				as desired. \qedhere
			\end{itemize}
		\end{proof}
		
		\begin{theorem}\label{thm:cartanpresentation}
			The Cartan subalgebra $\U_\A^0(\sl_2)$ has a presentation given by the generators $\qbinom{K;\hskip1pt c}{t}$ for $c\in \Z$, $t\in \N_0$, the multiplication relation in Theorem~\ref{thm:multiplication_rule}, and the relation in Proposition~\ref{prop:cplus2}.
		\end{theorem}
		\begin{proof}
			We need to show that the product of any two generators $\qbinom{K;\hskip1pt c}{t}\qbinom{K;\hskip1pt d}{s}$ can be written as an $\A$-linear combination of Lusztig's elements from the relation in Theorem~\ref{thm:multiplication_rule}, and the relation in Proposition~\ref{prop:cplus2}. Notice that by Theorem~\ref{thm:multiplication_rule} this is true if $t-c+b \ge 0$
			and $ s-b+c \ge 0$. Moreover $\qbinom{K;\hskip1pt c}{0} = 1$.
			Therefore, it suffices to show that one can express any $\qbinom{K;\hskip1pt c}{t}$ for $c \in \Z$ and $t \in \N$ as a linear combination of elements in the $\A$-basis $\B$ in Proposition~\ref{prop:binomialbasis}(1).
			
			We do this by induction on $t \in \N$. We show that the base case $t=1$ follows by a second induction on $c$, where the base cases $c=0,1$ hold by definition. If $c\ge 2$, assume that the result is true for $c'<c$. Then we can use Proposition~\ref{prop:cplus2} to write 
			$\qbinom{K;\hskip1pt c}{1}$ as an $\A$-linear combination of $\qbinom{K;c-1}{1}$ and $\qbinom{K;\hskip1pt c-2}{1}$, and so we are done by induction hypothesis on $c$. If $c<0$, we may assume that the result is true for $c'>c$. We then rearrange the relation in Proposition~\ref{prop:cplus2} to write 
			\[\qbinom{K;c}{1}=(q+q^{-1})\qbinom{K;c+1}{1}-\qbinom{K;c+2}{1}\]
			so again by induction hypothesis on $c$, we are done.
			
			Let us now consider the inductive step for $t>1$. By the induction hypothesis on $t$, the last term $\qbinom{K;\hskip1pt c}{t-2}$ in Proposition~\ref{prop:cplus2} can be expressed as an $\A$-linear combination of elements in $\B$. The relation without that term has the same form as for the case $t=1$, and the inductive argument on $c$ for the $t=1$ case carries over verbatim.
		\end{proof}

		\begin{corollary}\label{cor:changeOfBasis}
			Lusztig's integral form $\U_\A(\sl_2)$ has a presentation given by the monomials $E^{(n)}\qbinom{K;\hskip1pt c}{t}F^{(m)}$ for $n,m,t\ge 0$, $c\in\Z$, relations (\ref{eq:Erelation}), (\ref{eq:Frelation}), (\ref{eq:EKrelation}), (\ref{eq:FKrelation}) and (\ref{eq:EFrelation}), together with Proposition~\ref{prop:cplus2} and the multiplication relation in Theorem~\ref{thm:multiplication_rule}.
		\end{corollary}
		\begin{proof}
			Relation (\ref{eq:Erelation}) (resp.~(\ref{eq:Frelation})) allows us to rewrite products of divided powers $E^{(n)}E^{(m)}$ (resp. $F^{(n)}F^{(m)}$) as $\A$-multiples of a single divided power. Theorem~\ref{thm:cartanpresentation} allows us to rewrite products of Lusztig's elements $\qbinom{K;\hskip1pt c}{t}$ as $\A$-linear combinations of other Lusztig's elements. Finally, relations (\ref{eq:EKrelation}), (\ref{eq:FKrelation}) and (\ref{eq:EFrelation}) allow us to write any product of monomials $E^{(n)}\qbinom{K;c}{t}F^{(m)}$ in terms of other such monomials.
		\end{proof}

		Even though we have stated our results for $\U_\A(\sl_2)$, these bases generalize to arbitrary (finite) rank. Consider Lusztig's form of the quantum group $\U_\A(\mathfrak{g})$ for an arbitrary Kac--Moody algebra $\mathfrak{g}$ with generalized Cartan matrix of finite rank $n$, that is, with Cartan datum $(I,\cdot)$ such that $|I|=n$. 
		(The bilinear form $\cdot$ plays no role in the following.) By 1.4.7 and 3.1.13 in~\cite{Lus10}, its Cartan subalgebra is isomorphic to $U^0_\A(\sl_2)^{\otimes n}$, and we obtain the following corollary.
		\begin{corollary}\label{cor:allranks}
			Let $\mathfrak{g}$ be a Kac--Moody algebra of finite rank $n$, and let $K_i^{\pm1}$ for $i=1,\ldots,n$ be the generators of the Cartan subalgebra of the quantum group $\U_q(\mathfrak{g})$. Then Lusztig's form for the Cartan subalgebra $U^0_\A(\mathfrak{g})$ has an $\A$-basis given by the elements
			\[\qbinom{K_1\hskip0.5pt ;c}{t_1}\qbinom{K_2\hskip0.5pt ;c_2}{t_2}\cdots \qbinom{K_r\hskip0.5pt ;c_r}{t_r}\]
			with $c_i\in \{0,\min(1,t_i)\}$. Furthermore, it has a presentation with these generators and the relations in Theorem~\ref{thm:cartanpresentation} for each set of elements $\qbinom{K_i\hskip0.5pt ;\hskip1pt c}{t}$, together with the commutativity relations $\qbinom{K_i\hskip0.5pt;\hskip1pt c}{t}\qbinom{K_j\hskip0.5pt;\hskip1pt b}{s}
			=\qbinom{K_j\hskip0.5pt;\hskip1pt b}{s}\qbinom{K_i\hskip0.5pt;\hskip1pt c}{t}$ for $1\le i < j \le n$.
		\end{corollary}
		
		\section*{Acknowledgements}
		\'Alvaro Guti\'errez was funded by a University of Bristol Research Training Support Grant. Micha\l~Szwej and Mark Wildon gratefuly acknowledge financial support from the Heilbronn Institute for Mathematical Research, Bristol, UK. We
		thank Christian Krattenthaler for an insightful conversation about hypergeometric series. The authors are very grateful to two anonymous referees for exceptionally helpful critical comments and suggestions.

	\end{document}